\documentclass[11pt]{article}

\usepackage{subfig}
\usepackage{epic}
\usepackage{amssymb}
\usepackage{latexsym}
\usepackage{tikz}
\usetikzlibrary{decorations.markings}

\usetikzlibrary{arrows}

\newtheorem{thm}{Theorem}
\newtheorem{ob}[thm]{Observation}

\newtheorem{lem}[thm]{Lemma}
\newtheorem{cor}[thm]{Corollary}

\newcommand{\cN}{{\cal N}}
\newcommand{\cC}{{\cal C}}
\newcommand{\cub}{{\rm cubic}}

\newcommand{\cp}{\, \Box \,}

\newcommand{\smallqed}{{\tiny ($\Box$)}}

\newcommand{\proof}{\noindent\textbf{Proof. }}

\newcommand{\qed}{$\Box$}
\newenvironment{unnumbered}[1]{\trivlist
\item [\hskip \labelsep {\bf #1}]\ignorespaces\it}{\endtrivlist}

\newcommand{\2}{ \vspace{0.2cm} }

\definecolor{DarkGreen}{rgb}{0.2, 0.6, 0.3}

\definecolor{electricindigo}{rgb}{0.44, 0.0, 1.0}

\let\oldenumerate\enumerate
\renewcommand{\enumerate}{
  \oldenumerate
  \setlength{\itemsep}{0.5pt}
  \setlength{\parskip}{0pt}
  \setlength{\parsep}{0pt}
}

\begin{document}

\title{Total Forcing and Zero Forcing in \\ Claw-Free Cubic Graphs}
\author{$^{1,2}$Randy Davila and $^{1}$Michael A. Henning\\
\\
$^1$Department of Pure and Applied Mathematics\\
University of Johannesburg \\
Auckland Park 2006, South Africa \\
\small {\tt Email: mahenning@uj.ac.za} \\
\\
$^2$Department of Mathematics and Statistics \\
University of Houston--Downtown \\
Houston, TX 77002, USA \\
\small {\tt Email: davilar@uhd.edu}
}

\date{}
\maketitle

\begin{abstract}
A dynamic coloring of the vertices of a graph $G$ starts with an initial subset $S$ of colored vertices, with all remaining vertices being non-colored. At each discrete time interval, a colored vertex with exactly one non-colored neighbor forces this non-colored neighbor to be colored. The initial set $S$ is called a forcing set (zero forcing set) of $G$ if, by iteratively applying the forcing process, every vertex in $G$ becomes colored. If the initial set $S$ has the added property that it induces a subgraph of $G$ without isolated vertices, then $S$ is called a total forcing set in $G$. The total forcing number of $G$, denoted $F_t(G)$, is the minimum cardinality of a total forcing set in $G$. We prove that if $G$ is a connected cubic graph of order~$n$ that has a spanning $2$-factor consisting of triangles, then $F_t(G) \le \frac{1}{2}n$. More generally, we prove that if $G$ is a connected, claw-free, cubic graph of order~$n \ge 6$, then $F_t(G) \le \frac{1}{2}n$, where a claw-free graph is a graph that does not contain $K_{1,3}$ as an induced subgraph. The graphs achieving equality in these bounds are characterized.
\end{abstract}

{\small \textbf{Keywords:} Zero forcing sets; total forcing sets; claw-free; cubic; cycle cover.} \\
\indent {\small \textbf{AMS subject classification: 05C69}}

\newpage
\section{Introduction}

A dynamic coloring of the vertices in a graph is a coloring of the vertex set which may change, or propagate, throughout the vertices during discrete time intervals. Of the dynamic colorings, the notion of \emph{forcing sets} (\emph{zero forcing sets}), and the associated graph invariant known as the \emph{forcing number} (\emph{zero forcing number}), are arguably the most prominent, see for example \cite{AIM-Workshop, Barioli13, zf_np, DaHe17+, Davila Kenter, girthTheorem, Genter1, Genter2, Edholm, Hogben16, LuTang16, zf_np2}.
In this paper, we continue the study of total forcing sets in graphs, where a \emph{total forcing set} is a forcing set that induces a subgraph without isolated vertices.

More formally, let $G$ be a graph with vertex set $V = V(G)$ and edge set $E = E(G)$. The \emph{forcing process} is defined in~\cite{DaHe17a,DaHe17b} as follows: Let $S \subseteq V$ be a set of initially ``colored" vertices, all other vertices are said to be ``non-colored". A vertex contained in $S$ is said to be $S$-colored, while a vertex not in $S$ is said to be $S$-uncolored. At each time step, if a colored vertex has exactly one non-colored neighbor, then this colored vertex \emph{forces} its non-colored neighbor to become colored. If $v$ is such a colored vertex, we say that $v$ is a \emph{forcing vertex}. We say that $S$ is a \emph{forcing set}, if by iteratively applying the forcing process, all of $V$ becomes colored. We call such a set $S$, an $S$-forcing set. In addition, if $S$ is an $S$-forcing set in G and $v$ is a $S$-colored vertex that forces a new vertex to be colored, then $v$ is an \emph{$S$-forcing vertex}. The cardinality of a minimum forcing set in $G$ is the \emph{forcing number} of $G$, denoted $F(G)$.

If $S$ is a forcing set in a graph $G$ and the subgraph $G[S]$ induced by $S$ contains no isolated vertex, then $S$ is a \emph{total forcing set}, abbreviated as a TF-\emph{set} of $G$. The \emph{total forcing number} of $G$, written $F_t(G)$, is the cardinality of a minimum TF-set in $G$. The concept of a total forcing set was first introduced and studied by Davila in~\cite{Davila} as a strengthening of the concept of zero forcing originally introduced by the AIM-minimum rank group in~\cite{AIM-Workshop}. In~\cite{DaHe17a}, the authors show that if $G$ is a connected graph of order~$n \ge 3$ with maximum degree~$\Delta$, then $F_t(G) \le ( \frac{\Delta}{\Delta +1} ) n$. If we restrict $G$ to be a tree, then this bound can be improved to $F_t(G) \le \frac{1}{\Delta}((\Delta - 1)n + 1)$, as shown in~\cite{DaHe17b}.

A graph is \emph{$F$-free} if it does not contain $F$ as an induced subgraph. In particular, if $F = K_{1,3}$, then the graph is \emph{claw-free}. An excellent survey of claw-free graphs has been written by Flandrin, Faudree, and Ryj\'{a}\v{c}ek~\cite{ffr}. Chudnovsky and Seymour recently attracted considerable interest in claw-free graphs due to their excellent series of papers in \textit{Journal of Combinatorial Theory} on this topic (see, for example, their paper~\cite{ChSeV}). Domination and total domination in claw-free, cubic graph has been extensively studied (see, for example,~\cite{DeHaHe17,fh04,FaHe08b,HeLo12,HeMa16,HeYe_book,Li11,SoHe10} and elsewhere).
In this paper, we study total forcing sets in cubic, claw-free graphs.

\subsection{Notation}

For notation and graph theory terminology, we in general follow~\cite{HeYe_book}. Specifically, let $G$ be a graph with vertex set $V(G)$ and edge set $E(G)$, and of order~$n(G) = |V(G)|$ and size $m(G) = |E(G)|$. A \emph{neighbor} of a vertex $v$ in $G$ is a vertex $u$ that is adjacent to $v$, that is, $uv\in E(G)$. The \emph{open neighborhood} of a vertex $v$ in $G$ is the set of neighbors of $v$, denoted $N_G(v)$.
The \emph{neighborhood} of a set $S$ of vertices of $G$ is the set $N_G(S) = \cup_{v \in S} N_G(v)$. We denote the \emph{degree} of $v$ in $G$ by $d_G(v)$. A \emph{cubic graph} (also called a $3$-\emph{regular graph}) is a graph in which every vertex has degree~$3$, while a \emph{cubic multigraph} is a multigraph in which every vertex has degree~$3$.

For a set of vertices $S \subseteq V(G)$, the subgraph induced by $S$ is denoted by $G[S]$. The subgraph obtained from $G$ by deleting all vertices in $S$ and all edges incident with vertices in $S$ is denoted by $G - S$.
We denote the path and cycle on $n$ vertices by $P_n$ and $C_n$, respectively. We denote by $K_n$ the \emph{complete graph} on $n$ vertices. A complete graph $K_3$ is called a \emph{triangle}. The complete graph on four vertices minus one edge is called a \emph{diamond}. A vertex is \emph{covered} by the cycle if it belongs to the cycle. A \emph{cycle cover} of a multigraph is a collection of vertex-disjoint cycles in the multigraph that cover all the vertices. A \emph{$2$-factor} of a graph $G$ is a spanning $2$-regular subgraph of $G$; that is, a $2$-factor of $G$ is a collection of cycles that contain all vertices of $G$. We use the standard notation $[k] = \{1,2,\ldots,k\}$.

\subsection{Diamond-Necklace}
\label{S:neck}

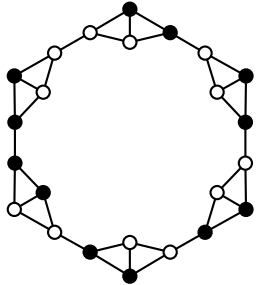
\begin{figure}[htb]
\tikzstyle{every node}=[circle, draw,fill=black!0, inner sep=0pt,minimum width=.175cm]
\begin{center}
\begin{tikzpicture}[thick,scale=.8]
\draw(0,0) { 
    +(0.00,3.33) -- +(0.48,3.06)
    +(0.00,3.33) -- +(0.01,2.56)
    +(0.01,2.56) -- +(0.48,3.06)
    +(0.48,3.06) -- +(0.67,3.71)
    +(0.67,3.71) -- +(0.00,3.33)
    +(0.01,2.56) -- +(0.01,1.88)
    +(0.01,1.88) -- +(0.00,1.11)
    +(0.00,1.11) -- +(0.48,1.39)
    +(0.01,1.88) -- +(0.48,1.39)
    +(0.48,1.39) -- +(0.67,0.73)
    +(0.67,0.73) -- +(0.00,1.11)
    +(0.67,0.73) -- +(1.26,0.40)
    +(1.26,0.40) -- +(1.92,0.56)
    +(1.92,0.56) -- +(1.92,0.00)
    +(1.92,0.00) -- +(1.26,0.40)
    +(1.92,0.56) -- +(2.59,0.40)
    +(2.59,0.40) -- +(1.92,0.00)
    +(2.59,4.05) -- +(1.92,3.89)
    +(1.92,3.89) -- +(1.92,4.44)
    +(1.92,4.44) -- +(2.59,4.05)
    +(1.92,3.89) -- +(1.26,4.05)
    +(1.26,4.05) -- +(1.92,4.44)
    +(1.26,4.05) -- +(0.67,3.71)
    +(2.59,4.05) -- +(3.17,3.71)
    +(3.17,3.71) -- +(3.37,3.06)
    +(3.37,3.06) -- +(3.85,3.33)
    +(3.85,3.33) -- +(3.17,3.71)
    +(3.37,3.06) -- +(3.84,2.56)
    +(3.85,3.33) -- +(3.84,2.56)
    +(3.84,2.56) -- +(3.84,1.88)
    +(3.84,1.88) -- +(3.37,1.39)
    +(3.37,1.39) -- +(3.85,1.11)
    +(3.85,1.11) -- +(3.84,1.88)
    +(3.37,1.39) -- +(3.17,0.73)
    +(3.17,0.73) -- +(3.85,1.11)
    +(3.17,0.73) -- +(2.59,0.40)
    +(0.48,1.39)  node[fill=black]{}
    +(0.00,1.11) node{}
    +(0.01,1.88) node[fill=black]{}
    +(0.48,3.06) node{}
    +(0.67,0.73) node{}
    +(0.01,2.56) node[fill=black]{}
    +(0.00,3.33) node[fill=black]{}
    +(0.67,3.71) node{}
    +(1.26,4.05) node{}
    +(1.92,3.89) node{}
    +(1.92,4.44) node[fill=black]{}
    +(2.59,4.05) node[fill=black]{}
    +(2.59,0.40) node{}
    +(1.26,0.40) node[fill=black]{}
    +(1.92,0.00) node[fill=black]{}
    +(1.92,0.56) node{}
    +(3.17,3.71) node{}
    +(3.37,3.06) node{}
    +(3.85,3.33) node[fill=black]{}
    +(3.84,2.56) node[fill=black]{}
    +(3.84,1.88) node{}
    +(3.17,0.73) node[fill=black]{}
    +(3.85,1.11) node[fill=black]{}
    +(3.37,1.39) node{}
   };
\end{tikzpicture}
\end{center}
\vskip -0.6 cm \caption{A diamond-necklace $N_6$.} \label{f:N6}
\end{figure}

In this paper, we shall adopt the terminology and notation of a diamond-necklace, coined by Henning and L\"{o}wenstein in~\cite{HeLo12}. For completeness, we repeat their definition here. For $k \ge 2$ an integer, let $N_k$ be the connected cubic graph constructed as follows. Take $k$ disjoint copies $D_1, D_2, \ldots, D_{k}$ of a diamond, where $V(D_i) = \{a_i,b_i,c_i,d_i\}$ and where $a_ib_i$ is the missing edge in $D_i$. Let $N_k$ be obtained from the disjoint union of these $k$ diamonds by adding the edges $\{a_ib_{i+1} \mid i \in [k-1] \}$  and adding the edge $a_{k}b_1$. We call $N_k$ a \emph{diamond-necklace with $k$ diamonds}. Let $\cN_\cub = \{ N_k \mid k \ge 2\}$. A diamond-necklace, $N_6$, with six diamonds is illustrated in Figure~\ref{f:N6}, where the darkened vertices form a TF-set in the graph.

\section{Main Results}

In this paper, we prove the following two results. Our first result shows that the total forcing number of a connected, claw-free, cubic graph in which every unit is a triangle-unit is at most one-half its order. A proof of Theorem~\ref{thm1} is presented in Section~\ref{S:proof1}.

\begin{thm}
\label{thm1}
If $G$ is a connected cubic graph of order~$n$ that has a spanning $2$-factor consisting of triangles, then $F_t(G) \le \frac{n}{2}$, with equality if and only if $G$ is the prism $C_3 \cp K_2$.
\end{thm}

More generally, we show that if we exclude the exceptional graph $K_4$, then the total forcing number of a connected, claw-free, cubic graph is always at most one-half its order, and we characterize the extremal graphs achieving equality in this bound. A proof of Theorem~\ref{thm2} is presented in Section~\ref{S:proof2}.

\begin{thm}
\label{thm2}
If $G \ne K_4$ is a connected, claw-free, cubic graph of order $n$, then $F_t(G) \le \frac{1}{2}n$ with equality if and only if $G \in \cN_\cub$ or $G$ is the prism $C_3 \cp K_2$.
\end{thm}

As a consequence of Theorem~\ref{thm2}, we have the following upper bound on the forcing number of a connected, claw-free, cubic graph. A proof of Theorem~\ref{thm3} is presented in Section~\ref{S:proof3}.

\begin{thm}
\label{thm3}
If $G \ne K_4$ is a connected, claw-free, cubic graph of order $n$, then $F(G) \le \frac{1}{2}n$ with equality if and only if $G$ is the diamond-necklace $N_2$ or the prism $C_3 \cp K_2$.
\end{thm}

As an immediate consequence of Theorem~\ref{thm3}, the forcing number of a connected, claw-free, cubic graph of order at least~$10$ is strictly less than one-half its order.

\begin{cor}
\label{cor1}
If $G$ is a connected, claw-free, cubic graph of order~$n \ge 10$, then $F(G) < \frac{1}{2}n$.
\end{cor}

We close with the following relationship between the total forcing and forcing numbers of a connected, claw-free, cubic graph.

\begin{thm}
\label{thm4}
If $G$ is a connected, claw-free, cubic graph of order $n$, then
\[
\frac{F_t(G)}{F(G)} \le 2,
\]
and this bound is asymptotically best possible. 
\end{thm}

\section{Known Results and Preliminary Lemma}

We shall need the following result observed, for example, in~\cite{Davila,DaHe17a} and elsewhere. 

\begin{ob}
\label{ob1}
If $G$ is an isolate-free graph, then $F(G) \le F_t(G) \le 2F(G)$.
\end{ob}

The following property of connected, claw-free, cubic graphs is established in~\cite{HeLo12}.

\begin{lem}{\rm (\cite{HeLo12})}
If $G \ne K_4$ is a connected, claw-free, cubic graph of order $n$, then the vertex set $V(G)$ can be uniquely partitioned into sets each of which induces a triangle or a diamond in $G$.
 \label{l:known}
\end{lem}

By Lemma~\ref{l:known}, the vertex set $V(G)$ of connected, claw-free, cubic graph $G \ne K_4$  can be uniquely partitioned into sets each of which induce a triangle or a diamond in $G$. Following the notation introduced in~\cite{HeLo12}, we refer to such a partition as a \emph{triangle}-\emph{diamond partition} of $G$, abbreviated $\Delta$-D-partition. We call every triangle and diamond induced by a set in our $\Delta$-D-partition a \emph{unit} of the partition. A unit that is a triangle is called a \emph{triangle-unit} and a unit that is a diamond is called a \emph{diamond-unit}. (We note that a triangle-unit is a triangle that does not belong to a diamond.) We say that two units in the $\Delta$-D-partition are \emph{adjacent} if there is an edge joining a vertex in one unit to a vertex in the other unit.
Before presenting proofs of our main results, we prove first that every graph in the family $\cN_\cub$ has total forcing number exactly one-half the order of the graph and has forcing number exactly one-fourth the order plus two. We remark that the following lemma may be of interest to those who study forcing as it relates to the linear algebraic minimum rank problem.

\begin{lem}
\label{lem1}
If $G \in \cN_\cub$ has order~$n$, then $F_t(G) = \frac{1}{2}n$ and $F(G) = \frac{1}{4}n + 2$.
%
\end{lem}
\proof Let $G \in \cN_\cub$ have order~$n$. Thus, $G$ is a diamond-necklace with $k$ diamonds for some $k \ge 2$, where $n = 4k$. We follow the notation introduced in Section~\ref{S:neck}. Hence, $G = N_k$ consists of $k$ vertex disjoint diamonds $D_1, D_2, \ldots, D_{k}$, where $V(D_i) = \{a_i,b_i,c_i,d_i\}$ and where $a_ib_i$ is the missing edge in $D_i$. Further, $G$ is obtained from these $k$ diamonds by adding the edges $\{a_ib_{i+1} \mid i \in [k-1] \}$ and adding the edge $a_{k}b_1$. Let $A = \{a_1,a_2,\ldots,a_k\}$ and $C = \{c_1,c_2,\ldots,c_k\}$.

(a) We first prove that $F_t(G) = 2k$. Let $S = (A \setminus \{a_1\}) \cup C \cup \{b_1\}$. (In the special case when $n = 6$, the set $S$ is illustrated by the darkened vertices in Figure~\ref{f:N6}, where here $G = N_6$ and $n = 24$.) The set $S$ is a TF-set of $G$ since the sequence $x_1,x_2,\ldots,x_{2k}$ of played vertices in the forcing process result in all vertices of $G$ colored, where $x_i$ denotes the forcing vertex played in the $i$th step of the process and where $x_1 = b_1$, $x_2 = d_1$, and $x_{2i+1} = a_{i+1}$ and $x_{2i+2} = d_{i+1}$ for $i \in [k-1]$; that is, the sequence of played vertices is given by $b_1, d_1, a_2, d_2, \ldots, a_k, b_k$. Thus, $F_t(G) \le |S| = 2k$.

We show next that $F_t(G) \ge 2k$. Let $S'$ be an arbitrary TF-set of $G$. We show that every diamond in the diamond-necklace contains at least two vertices of $S'$. If neither $c_i$ nor $d_i$ belong to the set $S'$ for some $i \in [k]$, then neither $c_i$ nor $d_i$ can be colored in the forcing process. Hence, at least one of $c_i$ and $d_i$ belongs to $S'$. Renaming vertices if necessary, we may assume that $c_i \in S'$. Since $G[S']$ is isolate-free, at least one neighbor of $c_i$ belongs to $S'$, implying that $|S' \cap V(D_i)| \ge 2$. Thus,
\[
|S'| = \sum_{i=1}^k |S' \cap V(D_i)| \ge 2k.
\]
\indent
Since $S'$ is an arbitrary TF-set of $G$, this implies that $F_t(G) \ge 2k$. As observed earlier, $F_t(G) \le 2k$. Consequently, $F_t(G) = 2k = \frac{1}{2}n$.

(b) We show next that $F(G) = k + 2$. In this case, we consider the set $S = C \cup \{b_1,a_2\}$. The set $S$ is a forcing set of $G$ since the sequence $x_1,x_2,\ldots,x_{3k+1}$ of played vertices in the forcing process result in all vertices of $G$ colored, where $x_i$ denotes the forcing vertex played in the $i$th step of the process and where $x_{3i-2} = a_{i+1}$, $x_{3i-1} = d_{i+1}$, and $x_{3i} = b_{i+1}$ for $i \in [k-1]$, and where $x_{3k-2} = a_1$; that is, the sequence of played vertices is given by $a_2,d_2,b_2, a_3,b_3,d_3, \ldots, a_k, d_k, b_k, a_1$. Thus, $F_t(G) \le |S| = k + 2$.

We show next that $F(G) \ge k + 2$. Let $S'$ be an arbitrary forcing set of $G$. At least one of $c_i$ and $d_i$ belong to $S'$ for every $i \in [k]$, and so every diamond in the diamond-necklace contains at least one vertex of $S'$. Renaming vertices if necessary, we may assume that $C \subseteq S'$. Since $S'$ is a forcing set, the first vertex $x_1$ played in the forcing process belongs to $S'$ and at least two neighbors of $x_1$ belong to $S'$. This implies that if $x_1 \in \{c_i,d_i\}$ for some $i \in [k]$, then $S'$ contains at least three vertices from the diamond $D_i$. Moreover, if $x_1 \in \{a_i,b_i\}$ for some $i \in [k]$, say $x_1 = b_i$ by symmetry, then $S'$ contains at least three vertices from the diamond $D_i$ or $S'$ contains at least two vertices from both diamonds $D_i$ and $D_{i+1}$ (where addition is taken modulo~$k$). Thus, one diamond in the diamond-necklace contains at least three vertices of $S'$ or two diamonds in the diamond-necklace both contain at least two vertices of $S'$, implying that $|S'| \ge k + 2$. Since $S'$ is an arbitrary TF-set of $G$, this implies that $F(G) \ge k+2$. As observed earlier, $F(G) \le k+2$. Consequently, $F(G) = k+2 = \frac{1}{4}n + 2$.~\qed

\section{Proof of Main Results}

\subsection{Proof of Theorem~\ref{thm1}}
\label{S:proof1}

In this section, we present a proof of Theorem~\ref{thm1}. Recall its statement.

\noindent \textbf{Theorem~\ref{thm1}}. \emph{If $G$ is a connected cubic graph of order~$n$ that has a spanning $2$-factor consisting of triangles, then $F_t(G) \le \frac{n}{2}$, with equality if and only if $G$ is the prism $C_3 \cp K_2$.
}

\medskip
\proof  Let $G$ be a connected cubic graph of order~$n$ that has a spanning $2$-factor consisting of triangles. We note that $G$ is a claw-free, cubic graph in which every unit is a triangle-unit. In particular, $G$ contains no diamond as a subgraph. We define the contraction multigraph of $G$, denoted $M_G$, to be the multigraph whose vertices correspond to the triangle-units in $G$ and where two vertices in $M_G$ are joined by the number of edges joining the corresponding triangle-units in $G$. By construction, $M_G$ has no loops but does possibly contain multiple edges. We note the order of $M_G$ is precisely the number of triangle-units in $G$. Since $G$ contains at least two triangle-units, $M_G$ has order at least~$2$. Every vertex in $M_G$ has degree~$3$, and so $M_G$ is a cubic multigraph. For each vertex $v$ in $M_G$, let $T_v$ denote the triangle-unit in $G$ associated with the vertex $v$.
Among all collections of vertex-disjoint cycles in $M_G$ (where we allow $2$-cycles), let $\cC$ be chosen so that the following holds. \2
\\
\indent (1) The number of vertices of $M_G$ covered by cycles that belong to $\cC$ is maximized. \\
\indent (2) Subject to (1), the number of cycles in $\cC$ is maximized.

\noindent
We now define a sequence of sets $S_1,S_2, \ldots$ as follows. Let $S_1$ be the set of all vertices covered by cycles in $\cC$; that is,
\[
S_1 = \bigcup_{C \in \cC} V(C).
\]

Since every vertex in the multigraph $M_G$ has degree~$3$, there is at least one cycle in $M_G$, implying that $S_1 \ne \emptyset$. Let $S_2$ be the set of all vertices in $M_G$ that do not belong to the set $S_1$, and are joined to $S_1$ by at least two edges (more precisely, by two or three edges). Let $S_3$ be the set of all vertices of $M_G$ not in $S_1 \cup S_2$ that are joined to $S_1 \cup S_2$ by at least two edges. More generally, if the set $S_i$ is defined for some $i \ge 1$ and
\[
S_{\le i} = \bigcup_{j=1}^i S_j,
\]
then let $S_{i+1}$ be the set of all vertices of $M_G$ not in $S_{\le i}$ that are joined to $S_{\le i}$ by at least two edges. Suppose that $V(M_G) \setminus S_{\le i} \ne \emptyset$ and $S_{i+1} = \emptyset$ for some $i \ge 0$; that is, suppose the set of vertices of $M_G$ that do not belong to the set $S_{\le i}$ is nonempty, but each such vertex is joined to $S_{\le i}$ by at most one edge. Let $M_{i+1}$ be the multigraph induced by the set of vertices of $M_G$ not in $S_{\le i}$. Since there is at most edge of $M_G$ that joins each vertex in $M_{i+1}$ to $S_{\le i}$, every vertex in the multigraph $M_{i+1}$ has degree at least~$2$, implying that there is a cycle in $M_{i+1}$. Adding this cycle to the collection of cycles in $\cC$ contradicts the maximality of $\cC$ (see condition~(2) of our choice of $\cC$). Hence, if $V(M_G) \setminus S_{\le i} \ne \emptyset$ for some $i \ge 0$, then $S_{i+1} \ne \emptyset$. Since the graph $G$ is finite, and so, $M_G$ is also finite, we note that for some integer $k \ge 0$, the following holds: $V(M_G) = S_{\le k}$ and $S_k \ne \emptyset$.

We now construct a TF-set of $G$ as follows. Let $C \colon v_1v_2\ldots v_\ell v_1$ be an arbitrary cycle in $\cC$, where $\ell \ge 2$. Let $T_1, T_2, \ldots, T_\ell$ be the triangle-units in $G$ associated with the vertices $v_1, v_2, \ldots, v_\ell$, respectively. Further, let $V(T_i) = \{v_{i1},v_{i2},v_{i3}\}$, where $v_{i2}v_{i+1,1}$ is an edge joining the triangle-units $T_i$ and $T_{i+1}$ in $G$ for $i \in [\ell]$, where addition is taken modulo~$\ell$. Thus, $C' \colon v_{11}v_{12}v_{21}v_{22} \ldots v_{\ell 1}v_{\ell 2} v_{11}$ is a cycle in $G$ that gives rise to the cycle $C$ in $M_G$. We remark that the cycle $C'$ may not be the unique such cycle. However, we select one such cycle $C'$ and fix this cycle. Relative to our choice of the cycle $C'$, we call the vertex $v_{i3}$ the \emph{free vertex} of the triangle-unit $T_i$. Let $n_C$ be the number of vertices that belong to triangle-units of $G$ that are associated with the cycle $C$ of $M_G$, and so, $n_C = 3\ell$. Next we color a subset $S$ of the vertices that belong to triangle-units of $G$ associated with the cycle $C$ of $M_G$ as follows.

\begin{figure}[htb]
\begin{center}
\begin{tikzpicture}[scale=.8,style=thick,x=1cm,y=1cm]
\def\vr{2.5pt} 
\path (-0.1,0.5) coordinate (v1);
\path (-0.1,-0.2) coordinate (u1);
\path (0,0.5) coordinate (v11);
\path (1.5,0.5) coordinate (v12);
\path (0.75,1.5) coordinate (v13);
\path (0.75,2.1) coordinate (v1p);
\path (3,0.5) coordinate (v21);
\path (4.5,0.5) coordinate (v22);
\path (3.75,1.5) coordinate (v23);
\path (3.75,2.1) coordinate (v2p);
\path (6,0.5) coordinate (v31);
\path (7.5,0.5) coordinate (v32);
\path (6.75,1.5) coordinate (v33);
\path (6.75,2.1) coordinate (v3p);
\path (9,0.5) coordinate (v41);
\path (10.5,0.5) coordinate (v42);
\path (9.75,1.5) coordinate (v43);
\path (9.75,2.1) coordinate (v4p);
\path (12,0.5) coordinate (v51);
\path (13.5,0.5) coordinate (v52);
\path (12.75,1.5) coordinate (v53);
\path (12.75,2.1) coordinate (v5p);
\path (15,0.5) coordinate (v61);
\path (16.5,0.5) coordinate (v62);
\path (15.75,1.5) coordinate (v63);
\path (15.75,2.1) coordinate (v6p);
\path (16.6,0.5) coordinate (v6);
\path (16.6,-0.2) coordinate (u6);
\draw (v11) -- (v12);
\draw (v11) -- (v13);
\draw (v12) -- (v13);
\draw (v21) -- (v22);
\draw (v21) -- (v23);
\draw (v22) -- (v23);
\draw (v31) -- (v32);
\draw (v31) -- (v33);
\draw (v32) -- (v33);
\draw (v41) -- (v42);
\draw (v41) -- (v43);
\draw (v42) -- (v43);
\draw (v51) -- (v52);
\draw (v51) -- (v53);
\draw (v52) -- (v53);
\draw (v61) -- (v62);
\draw (v61) -- (v63);
\draw (v62) -- (v63);
\draw (v12) -- (v21);
\draw (v22) -- (v31);
\draw (v32) -- (v41);
\draw (v42) -- (v51);
\draw (v52) -- (v61);
\draw (v62) -- (v11);
\draw (v13) -- (v1p);
\draw (v23) -- (v2p);
\draw (v33) -- (v3p);
\draw (v43) -- (v4p);
\draw (v53) -- (v5p);
\draw (v63) -- (v6p);
\draw (u1) -- (u6);
\draw (v11) [fill=black] circle (\vr);
\draw (v12) [fill=black] circle (\vr);
\draw (v13) [fill=black] circle (\vr);
\draw (v21) [fill=white] circle (\vr);
\draw (v22) [fill=white] circle (\vr);
\draw (v23) [fill=white] circle (\vr);
\draw (v31) [fill=black] circle (\vr);
\draw (v32) [fill=black] circle (\vr);
\draw (v33) [fill=black] circle (\vr);
\draw (v41) [fill=white] circle (\vr);
\draw (v42) [fill=white] circle (\vr);
\draw (v43) [fill=white] circle (\vr);
\draw (v51) [fill=black] circle (\vr);
\draw (v52) [fill=black] circle (\vr);
\draw (v53) [fill=black] circle (\vr);
\draw (v61) [fill=white] circle (\vr);
\draw (v62) [fill=white] circle (\vr);
\draw (v63) [fill=white] circle (\vr);
%
\draw (v1) to[out=180,in=180, distance=1cm] (u1);
\draw (v6) to[out=0,in=0, distance=1cm] (u6);
%
\draw (0,0.25) node {{\small $v_{11}$}};
\draw (1.5,0.25) node {{\small $v_{12}$}};
\draw (1.2,1.5) node {{\small $v_{13}$}};
\draw (3,0.25) node {{\small $v_{21}$}};
\draw (4.5,0.25) node {{\small $v_{22}$}};
\draw (4.2,1.5) node {{\small $v_{23}$}};
\draw (6,0.25) node {{\small $v_{31}$}};
\draw (7.5,0.25) node {{\small $v_{32}$}};
\draw (7.2,1.5) node {{\small $v_{33}$}};
\draw (9,0.25) node {{\small $v_{41}$}};
\draw (10.5,0.25) node {{\small $v_{42}$}};
\draw (10.2,1.5) node {{\small $v_{43}$}};
\draw (12,0.25) node {{\small $v_{51}$}};
\draw (13.5,0.25) node {{\small $v_{52}$}};
\draw (13.2,1.5) node {{\small $v_{53}$}};
\draw (15,0.25) node {{\small $v_{61}$}};
%
\draw (16.5,0.25) node {{\small $v_{62}$}};
\draw (16.2,1.5) node {{\small $v_{63}$}};
\draw (0.75,0.9) node {{\small $T_1$}};
\draw (3.75,0.9) node {{\small $T_2$}};
\draw (6.75,0.9) node {{\small $T_3$}};
\draw (9.75,0.9) node {{\small $T_4$}};
\draw (12.75,0.9) node {{\small $T_5$}};
\draw (15.75,0.9) node {{\small $T_6$}};
\end{tikzpicture}
\end{center}
\vskip -0.25cm
\caption{The coloring of the triangle-units when $\ell = 6$.} \label{f:Ceven}
\end{figure}
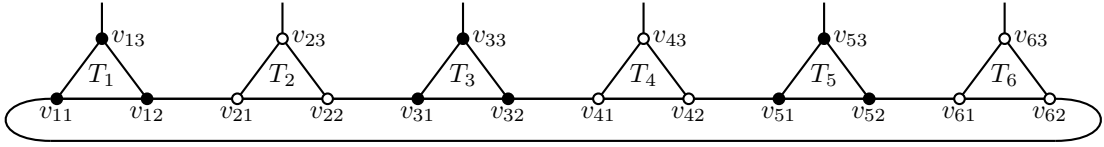

If $\ell \ge 2$ is even, then we color all vertices that belong to the triangle-units $T_{2i-1}$ where $i \in [\ell/2]$. Let $S_C$ denote the resulting set of colored vertices of $G$ associated with the cycle~$C$. (We illustrate the case when $\ell = 6$ in Figure~\ref{f:Ceven}, where the vertices in $S_C$ are darkened.) We note that exactly one-half the vertices in all triangle-units of $G$ associated with the cycle $C$ of $M_G$ are colored; that is, $|S_C| = n_C/2$.

\begin{figure}[htb]
\begin{center}
\begin{tikzpicture}[scale=.8,style=thick,x=1cm,y=1cm]
\def\vr{2.5pt} 
\path (-0.1,0.5) coordinate (v1);
\path (-0.1,-0.2) coordinate (u1);
\path (0,0.5) coordinate (v11);
\path (1.5,0.5) coordinate (v12);
\path (0.75,1.5) coordinate (v13);
\path (0.75,2.1) coordinate (v1p);
\path (3,0.5) coordinate (v21);
\path (4.5,0.5) coordinate (v22);
\path (3.75,1.5) coordinate (v23);
\path (3.75,2.1) coordinate (v2p);
\path (6,0.5) coordinate (v31);
\path (7.5,0.5) coordinate (v32);
\path (6.75,1.5) coordinate (v33);
\path (6.75,2.1) coordinate (v3p);
\path (9,0.5) coordinate (v41);
\path (10.5,0.5) coordinate (v42);
\path (9.75,1.5) coordinate (v43);
\path (9.75,2.1) coordinate (v4p);
\path (12,0.5) coordinate (v51);
\path (13.5,0.5) coordinate (v52);
\path (12.75,1.5) coordinate (v53);
\path (12.75,2.1) coordinate (v5p);
\path (13.6,0.5) coordinate (v5);
\path (13.6,-0.2) coordinate (u5);
\draw (v11) -- (v12);
\draw (v11) -- (v13);
\draw (v12) -- (v13);
\draw (v21) -- (v22);
\draw (v21) -- (v23);
\draw (v22) -- (v23);
\draw (v31) -- (v32);
\draw (v31) -- (v33);
\draw (v32) -- (v33);
\draw (v41) -- (v42);
\draw (v41) -- (v43);
\draw (v42) -- (v43);
\draw (v51) -- (v52);
\draw (v51) -- (v53);
\draw (v52) -- (v53);
\draw (v12) -- (v21);
\draw (v22) -- (v31);
\draw (v32) -- (v41);
\draw (v42) -- (v51);
\draw (v52) -- (v11);
\draw (v13) -- (v1p);
\draw (v23) -- (v2p);
\draw (v33) -- (v3p);
\draw (v43) -- (v4p);
\draw (v53) -- (v5p);
\draw (u1) -- (u5);
\draw (v11) [fill=white] circle (\vr);
\draw (v12) [fill=black] circle (\vr);
\draw (v13) [fill=black] circle (\vr);
\draw (v21) [fill=black] circle (\vr);
\draw (v22) [fill=white] circle (\vr);
\draw (v23) [fill=black] circle (\vr);
\draw (v31) [fill=white] circle (\vr);
\draw (v32) [fill=white] circle (\vr);
\draw (v33) [fill=white] circle (\vr);
\draw (v41) [fill=black] circle (\vr);
\draw (v42) [fill=black] circle (\vr);
\draw (v43) [fill=black] circle (\vr);
\draw (v51) [fill=white] circle (\vr);
\draw (v52) [fill=white] circle (\vr);
\draw (v53) [fill=white] circle (\vr);
%
\draw (v1) to[out=180,in=180, distance=1cm] (u1);
\draw (v5) to[out=0,in=0, distance=1cm] (u5);
%
\draw (0,0.25) node {{\small $v_{11}$}};
\draw (1.5,0.25) node {{\small $v_{12}$}};
\draw (1.2,1.5) node {{\small $v_{13}$}};
\draw (3,0.25) node {{\small $v_{21}$}};
\draw (4.5,0.25) node {{\small $v_{22}$}};
\draw (4.2,1.5) node {{\small $v_{23}$}};
\draw (6,0.25) node {{\small $v_{31}$}};
\draw (7.5,0.25) node {{\small $v_{32}$}};
\draw (7.2,1.5) node {{\small $v_{33}$}};
\draw (9,0.25) node {{\small $v_{41}$}};
\draw (10.5,0.25) node {{\small $v_{42}$}};
\draw (10.2,1.5) node {{\small $v_{43}$}};
\draw (12,0.25) node {{\small $v_{51}$}};
\draw (13.5,0.25) node {{\small $v_{52}$}};
\draw (13.2,1.5) node {{\small $v_{53}$}};
\draw (0.75,0.9) node {{\small $T_1$}};
\draw (3.75,0.9) node {{\small $T_2$}};
\draw (6.75,0.9) node {{\small $T_3$}};
\draw (9.75,0.9) node {{\small $T_4$}};
\draw (12.75,0.9) node {{\small $T_5$}};
\end{tikzpicture}
\end{center}
\vskip -0.25cm
\caption{The coloring of the triangle-units when $\ell = 5$.} \label{f:Codd}
\end{figure}
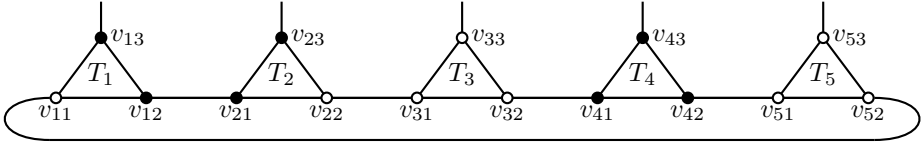

If $\ell \ge 3$ is odd, then we color the four vertices $v_{12},v_{13},v_{21},v_{23}$ and if $\ell \ge 5$ we further color all vertices that belong to the triangle-units $T_{2i}$ where $i \in  [(\ell-1)/2] \setminus \{1\}$. Let $S_C$ denote the resulting set of colored vertices of $G$ associated with the cycle~$C$. (We illustrate the case when $\ell = 5$ in Figure~\ref{f:Codd}, where the vertices in $S_C$ are darkened.) We note that $(n_C - 1)/2$ of the vertices in all triangle-units of $G$ associated with the cycle $C$ of $M_G$ are colored; that is, $|S_C| = (n_C - 1)/2$.

Let $S$ be the resulting set of all colored vertices associated with the cycles $C$ in $\cC$; that is,
\[
S = \bigcup_{C\in \cC} S_C
\hspace*{0.75cm} \mbox{and}  \hspace*{0.75cm}
|S| = \sum_{C\in \cC} |S_C| \le \frac{1}{2}\sum_{C\in \cC} n_C.
\]

We claim that $S$ is a TF-set in $G$. By construction, $G[S]$ contains no isolated vertex. Hence it suffices for us to prove that $S$ is a forcing set of $G$. For this purpose we define $k$ stages in the forcing process, where we recall that $k$ is the nonnegative integer such that $V(M_G) = S_{\le k}$ and $S_k \ne \emptyset$. For $i \in [k]$, let $V_i$ be the set of all vertices that belong to a triangle-unit of $G$ associated with a vertex of $M_G$ that belongs to the set $S_i$, and let \[
V_{\le i} = \bigcup_{j=1}^i V_j.
\]
We note that $V(G) = V_{\le k}$ and
\[
|V_1| = 3|S_1| = \sum_{C \in \cC} n_C \ge 2|S|,
\]
and so
\[
|S| \le \frac{1}{2}|V_1| \le \frac{1}{2}n.
\]

We first define Stage~1 of the forcing process. In this first stage we color all vertices in $V_1$, starting with the initial set $S$ of colored vertices, in the following way. We consider the cycles in $\cC$ sequentially (in any order) and color the vertices that belong to a triangle-unit of $G$ associated with such a cycle as follows. Let $C \colon v_1v_2\ldots v_\ell v_1$ in $\cC$, where $\ell \ge 2$, be an arbitrary cycle in $\cC$.

Suppose that $\ell \ge 2$ is even. Using our earlier notation, we note that the $S$-colored vertex $v_{2i-1,2}$ forces the vertex $v_{2i,1}$ to be colored for all $i \in [\ell/2]$. Further, the $S$-colored vertex $v_{11} \in S$ forces the vertex $v_{\ell,2}$ to be colored, while the $S$-colored vertex $v_{2i-1,1} \in S$ forces the vertex $v_{2i-2,2}$ to be colored for all $i \in [\ell/2] \setminus \{1\}$. In this way, all vertices $v_{ij}$, where $j \in [2]$, are colored in the forcing process. The newly colored vertex $v_{2i,1}$ now forces the vertex $v_{2i,3}$ to be colored for all $i \in [\ell/2]$. Thus if $\ell$ is even, then all vertices that belong to a triangle-unit of $G$ associated with the cycle $C$ of $M_G$ are colored starting with the initial set $S$ of colored vertices.

Suppose next that $\ell \ge 3$ is odd. Using our earlier notation, we note that the $S$-colored vertex $v_{12}$ forces the vertex $v_{11}$ to be colored, and the $S$-colored vertex $v_{21}$ forces the vertex $v_{22}$ to be colored. In this way, all six vertices in $T_1$ and $T_2$ are colored. The colored vertex $v_{2i,2}$ now forces the vertex $v_{2i+1,1}$ to be colored for all $i \in [(\ell - 1)/2]$. Further, the colored vertex $v_{11} \in S$ forces the vertex $v_{\ell,2}$ to be colored, while the $S$-colored vertex $v_{2i,1} \in S$ forces the vertex $v_{2i-1,2}$ to be colored for all $i \in [(\ell-1)/2] \setminus \{1\}$. In this way, all vertices $v_{ij}$, where $j \in [2]$, are colored in the forcing process. The newly colored vertex $v_{2i+1,1}$ now forces the vertex $v_{2i+1,3}$ to be colored for all $i \in [(\ell-1)/2]$. Thus if $\ell$ is odd, then all vertices that belong to a triangle-unit of $G$ associated with the cycle $C$ of $M_G$ are colored starting with the initial set $S$ of colored vertices. Thus, after Stage~1 of the forcing process all vertices in $V_1$ are colored.

We next define Stage~2 of the forcing process. In this second stage we color all vertices in $V_2$ as follows. Let $v$ be an arbitrary vertex in $V_2$. Thus, the vertex $v$ belongs to a triangle-unit, $T_v$ say, in $G$ that is associated with a vertex $v' \in S_2$ in $M_G$. By definition, the vertex $v'$ is joined to $S_1$ in $M_G$ by two or three edges. This implies that there are two or three edges joining $T_v$ to vertices in $V_1$. Let $a_2, b_2, c_2$ be the three vertices in the triangle-unit $T_v$ in $G$. We note that $v \in \{a_2,b_2,c_2\}$. Renaming vertices if necessary, we may assume that $a_2$ has a neighbor, $a_1$ say, in $V_1$ and $b_2$ has a neighbor, $b_1$ say, in $V_1$. Since every vertex of $V_1$ belongs to a triangle-unit in $G[V_1]$, we note that every vertex of $V_1$ has at most one neighbor not in $V_1$. In particular, $a_2$ is the only uncolored neighbor of $a_1$ and $b_2$ is the only uncolored neighbor of $b_1$ upon completion of Stage~1. Further, $a_1 \ne b_1$. The vertex $a_1$ now forces the vertex $a_2$ to be colored, and the vertex $b_1$ now forces the vertex $b_2$ to be colored. The newly colored vertex $a_2$ now forces the third vertex of $T_v$, namely $c_2$, to be colored. Thus, all three vertices in the triangle-unit $T_v$ in $G$ become colored in Stage~2. In particular, $v$ is colored in this stage of the forcing process. Since $v$ is an arbitrary vertex in $V_2$, after Stage~2 of the forcing process all vertices in $V_{\le 2} = V_1 \cup V_2$ are colored.

In general, suppose that Stage~$i$ in the forcing process has been defined for some $i$ where $1 \le i < k$, and that after Stage~$i$ all vertices in $V_{\le i}$ are colored. We now define the $(i+1)$th stage in the forcing process. In this stage we color all vertices in $V_{i+1}$ as follows. Let $v$ be an arbitrary vertex in $V_{i+1}$. Thus, the vertex $v$ belongs to a triangle-unit, $T_v$ say, in $G$ that is associated with a vertex $v' \in S_{i+1}$ in $M_G$. By definition, the vertex $v'$ is joined to $S_{\le i}$ in $M_G$ by two or three edges. This implies that there are two or three edges joining $T_v$ to vertices in $V_{\le i}$. Let $a_{i+1}, b_{i+1}, c_{i+1}$ be the three vertices in the triangle-unit $T_v$ in $G$. We note that $v \in \{a_{i+1},b_{i+1},c_{i+1}\}$. Renaming vertices if necessary, we may assume that $a_{i+1}$ has a neighbor, $a$ say, in $V_{\le i}$ and $b_{i+1}$ has a neighbor, $b$ say, in $V_{\le i}$.

Since every vertex of $V_{\le i}$ belongs to a triangle-unit in $G[V_{\le i}]$, we note that every vertex of $V_{\le i}$ has at most one neighbor not in $V_{\le i}$. In particular, $a_{i+1}$ is the only uncolored neighbor of $a$ and $b_{i+1}$ is the only uncolored neighbor of $b$ upon completion of Stage~$i$. Further, $a \ne b$. The vertex $a$ now forces the vertex $a_{i+1}$ to be colored, and the vertex $b$ now forces the vertex $b_{i+1}$ to be colored. The newly colored vertex $a_{i+1}$ now forces the third vertex of $T_v$, namely $c_{i+1}$, to be colored. Thus, all three vertices in the triangle-unit $T_v$ in $G$. In particular, $v$ is colored in this stage of the forcing process. Since $v$ is an arbitrary vertex in $V_{i+1}$, after Stage~$i+1$ of the forcing process all vertices in $V_{\le i+1}$ are colored. In particular, after Stage~$k$ of the forcing process all vertices in $V(G) = V_{\le k}$ are colored. Thus, the set $S$ is a TF-set of $G$, implying by our earlier observations that
\[
F_t(G) \le |S| \le \frac{1}{2}|V_1| \le \frac{1}{2}n.
\]

Suppose, next, that $F_t(G) = n/2$. Thus, we must have equality throughout the above inequality chain, implying in particular that $F_t(G) = |S| = |V_1|/2$ and $k = 1$. This in turn implies that $\cC$ is a cycle cover in $M_G$, and so every vertex of $M_G$ belongs to a cycle in $\cC$. Further, for every cycle $C$ in $\cC$, we have $|S_C| = n_C/2$, and so every cycle in $\cC$ has even length. Thus, the coloring of the triangle-units of the cycle $C$ is as illustrated in Figure~\ref{f:Ceven}.

Let $C \colon v_1v_2\ldots v_\ell v_1$ be an arbitrary cycle in $\cC$, where $\ell \ge 2$ is even. As before, let $T_1, T_2, \ldots, T_\ell$ be the triangle-units in $G$ associated with the vertices $v_1, v_2, \ldots, v_\ell$, respectively, and let $V(T_i) = \{v_{i1},v_{i2},v_{i3}\}$, where $C' \colon v_{11}v_{12}v_{21}v_{22} \ldots v_{\ell 1}v_{\ell 2} v_{11}$ is the chosen cycle in $G$ that gives rise to the cycle $C$ in $M_G$. Recall that the vertex $v_{i3}$ is called the free vertex of the triangle-unit $T_i$ for $i \in [\ell]$ relative to $C'$. We now consider the free vertex $v_{13}$. Let $v$ be the neighbor of $v_{13}$ that does not belong to the triangle-unit $T_1$. We note that $v$ is a free vertex associated with some cycle $C^*$ in $G$, where possibly $C' = C^*$. Let $T_v$ be the triangle-unit in $G$ that contains~$v$. We note that either all three vertices of $T_v$ are colored (that is, belong to $S$) or none are colored.

If $v$ is an $S$-colored vertex, then the set $S' = S \setminus \{v\}$ is a TF-set of $G$, since as the first vertex played in the forcing process starting with the set $S'$ we play the vertex $v_{13}$ which forces $v$ to be colored, and then proceed exactly with the forcing process starting with the set $S$ to color all remaining vertices, implying that $F_t(G) \le |S'| < |S| = n/2$, a contradiction. Therefore, $v$ is not an $S$-colored vertex, and so no vertex in $T_v$ is $S$-colored.

Suppose that $T_v$ is neither the triangle-unit $T_2$ nor the triangle-unit $T_\ell$ (where possibly $T_2 = T_\ell$). Let $x$ and $y$ be the two vertices of $T_v$ different from $v$, and let $x^*$ and $y^*$ be the neighbors of $x$ and $y$, respectively, that belong to the cycle $C^*$. Since $T_v \ne T_2$ and $T_v \ne T_\ell$, we note that neither $x^*$ nor $y^*$ belong to the triangle-units $T_2$ or $T_\ell$. We note, however, that all three vertices in the triangle-unit containing $x^*$ (respectively, $y^*$) are colored (that is, belong to $S$). In this case, the set $S' = S \setminus \{v_{13}\}$ is a TF-set of $G$, since as the first four vertices played in the forcing process starting with the set $S'$ we play first the vertex $x^*$ which forces $x$ to be colored, second we play the vertex  $y^*$ which forces $y$ to be colored, third we play the vertex  $x$ which forces $v$ to be colored, and fourth we play the vertex $v$ which forces $v_{13}$ to be colored, and then proceed with the forcing process starting with the set $S$ to color all remaining vertices, implying once again that $F_t(G) \le |S'| < |S| \le n/2$, a contradiction.

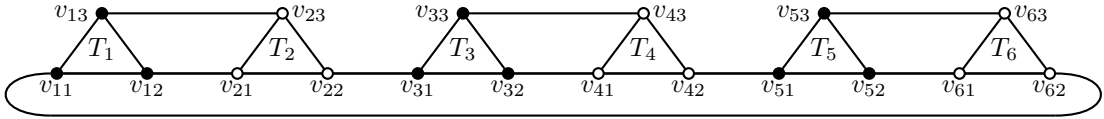
\begin{figure}[htb]
\begin{center}
\begin{tikzpicture}[scale=.8,style=thick,x=1cm,y=1cm]
\def\vr{2.5pt} 
\path (-0.1,0.5) coordinate (v1);
\path (-0.1,-0.2) coordinate (u1);
\path (0,0.5) coordinate (v11);
\path (1.5,0.5) coordinate (v12);
\path (0.75,1.5) coordinate (v13);
\path (3,0.5) coordinate (v21);
\path (4.5,0.5) coordinate (v22);
\path (3.75,1.5) coordinate (v23);
\path (6,0.5) coordinate (v31);
\path (7.5,0.5) coordinate (v32);
\path (6.75,1.5) coordinate (v33);
\path (9,0.5) coordinate (v41);
\path (10.5,0.5) coordinate (v42);
\path (9.75,1.5) coordinate (v43);
\path (12,0.5) coordinate (v51);
\path (13.5,0.5) coordinate (v52);
\path (12.75,1.5) coordinate (v53);
\path (15,0.5) coordinate (v61);
\path (16.5,0.5) coordinate (v62);
\path (15.75,1.5) coordinate (v63);
\path (16.6,0.5) coordinate (v6);
\path (16.6,-0.2) coordinate (u6);
\draw (v11) -- (v12);
\draw (v11) -- (v13);
\draw (v12) -- (v13);
\draw (v21) -- (v22);
\draw (v21) -- (v23);
\draw (v22) -- (v23);
\draw (v31) -- (v32);
\draw (v31) -- (v33);
\draw (v32) -- (v33);
\draw (v41) -- (v42);
\draw (v41) -- (v43);
\draw (v42) -- (v43);
\draw (v51) -- (v52);
\draw (v51) -- (v53);
\draw (v52) -- (v53);
\draw (v61) -- (v62);
\draw (v61) -- (v63);
\draw (v62) -- (v63);
\draw (v12) -- (v21);
\draw (v22) -- (v31);
\draw (v32) -- (v41);
\draw (v42) -- (v51);
\draw (v52) -- (v61);
\draw (v62) -- (v11);
%
\draw (v13) -- (v23);
\draw (v33) -- (v43);
\draw (v53) -- (v63);
\draw (u1) -- (u6);
\draw (v11) [fill=black] circle (\vr);
\draw (v12) [fill=black] circle (\vr);
\draw (v13) [fill=black] circle (\vr);
\draw (v21) [fill=white] circle (\vr);
\draw (v22) [fill=white] circle (\vr);
\draw (v23) [fill=white] circle (\vr);
\draw (v31) [fill=black] circle (\vr);
\draw (v32) [fill=black] circle (\vr);
\draw (v33) [fill=black] circle (\vr);
\draw (v41) [fill=white] circle (\vr);
\draw (v42) [fill=white] circle (\vr);
\draw (v43) [fill=white] circle (\vr);
\draw (v51) [fill=black] circle (\vr);
\draw (v52) [fill=black] circle (\vr);
\draw (v53) [fill=black] circle (\vr);
\draw (v61) [fill=white] circle (\vr);
\draw (v62) [fill=white] circle (\vr);
\draw (v63) [fill=white] circle (\vr);
%
\draw (v1) to[out=180,in=180, distance=1cm] (u1);
\draw (v6) to[out=0,in=0, distance=1cm] (u6);
%
\draw (0,0.25) node {{\small $v_{11}$}};
\draw (1.5,0.25) node {{\small $v_{12}$}};
\draw (0.25,1.5) node {{\small $v_{13}$}};
\draw (3,0.25) node {{\small $v_{21}$}};
\draw (4.5,0.25) node {{\small $v_{22}$}};
\draw (4.2,1.5) node {{\small $v_{23}$}};
\draw (6,0.25) node {{\small $v_{31}$}};
\draw (7.5,0.25) node {{\small $v_{32}$}};
\draw (6.25,1.5) node {{\small $v_{33}$}};
\draw (9,0.25) node {{\small $v_{41}$}};
\draw (10.5,0.25) node {{\small $v_{42}$}};
\draw (10.2,1.5) node {{\small $v_{43}$}};
\draw (12,0.25) node {{\small $v_{51}$}};
\draw (13.5,0.25) node {{\small $v_{52}$}};
\draw (12.25,1.5) node {{\small $v_{53}$}};
\draw (15,0.25) node {{\small $v_{61}$}};
%
\draw (16.5,0.25) node {{\small $v_{62}$}};
\draw (16.2,1.5) node {{\small $v_{63}$}};
\draw (0.75,0.9) node {{\small $T_1$}};
\draw (3.75,0.9) node {{\small $T_2$}};
\draw (6.75,0.9) node {{\small $T_3$}};
\draw (9.75,0.9) node {{\small $T_4$}};
\draw (12.75,0.9) node {{\small $T_5$}};
\draw (15.75,0.9) node {{\small $T_6$}};
\end{tikzpicture}
\end{center}
\vskip -0.25cm
\caption{The graph $G$ when $\ell = 6$.} \label{f:G}
\end{figure}

Therefore, $T_v$ is the triangle-unit $T_2$ or the triangle-unit $T_\ell$. Since $G$ is connected and the vertex $v_{13}$ is an arbitrary free vertex in $G$, this implies that $\cC$ consists of exactly one cycle, that is, $\cC = \{C\}$ and $|\cC| = 1$. By symmetry, we may assume that $T_v = T_2$. Suppose that $C$ has length $\ell \ge 4$. We now consider the free vertex $v_{33}$ that belongs to the triangle-unit $T_3$. Analogously as with the vertex $v_{13}$, we show that the free vertex $v_{33}$ is adjacent with the free vertex $v_{43}$ that belongs to the triangle-unit $T_4$. More generally, we show that the free vertex $v_{2i+1,3}$ that belongs to the triangle-unit $T_{2i+1}$ is adjacent with the free vertex $v_{2i+2,3}$ that belongs to the triangle-unit $T_{2i+2}$ for all $i \in [\ell - 1]$. Thus the graph $G$ is determined. We illustrate the graph $G$ in the special case when $\ell = 6$ in Figure~\ref{f:G}, where the vertices in $S$ are darkened.

\medskip

However, we can now replace the cycle $C \colon v_1v_2\ldots v_\ell v_1$ in $\cC$ that belongs to the multigraph $M_G$ with the $2$-cycles $v_1v_2v_1$, $v_3v_4v_3, \ldots, v_{\ell - 1} v_\ell v_{\ell - 1}$ of $M_G$. Thus we could have chosen the cycle cover $\cC$ to consist of $\ell/2 \ge 2$ vertex-disjoint cycles, contradicting condition~(2) of our choice of $\cC$. Therefore, $\ell = 2$. Thus, $|\cC| = 1$ and the cycle in $\cC$ is a $2$-cycle, implying that $G$ contains exactly two triangle-units and is therefore the prism $C_3 \cp K_2$. Therefore if $F_t(G) = n/2$, then we have shown that $G = C_3 \cp K_2$.~\qed

\subsection{Proof of Theorem~\ref{thm2}}
\label{S:proof2}

In this section, we present a proof of Theorem~\ref{thm2}. Recall its statement.

\noindent \textbf{Theorem~\ref{thm2}}. \emph{If $G \ne K_4$ is a connected, claw-free, cubic graph of order $n$, then $F_t(G) \le \frac{1}{2}n$ with equality if and only if $G \in \cN_\cub$ or $G$ is the prism $C_3 \cp K_2$.
}

\begin{figure}[htb]
\begin{center}
\begin{tikzpicture}[scale=.8,style=thick,x=1cm,y=1cm]
\def\vr{2.75pt} 
\path (0,0) coordinate (b);
\path (0,1.5) coordinate (a);
\path (1,0.75) coordinate (c);
\path (2,0.75) coordinate (d);
\path (3,0) coordinate (f);
\path (3,1.5) coordinate (e);
\path (5.25,0.75) coordinate (g);
\path (6.75,0.75) coordinate (j);
\path (8.25,0.75) coordinate (k);
\path (9.75,0.75) coordinate (n);
\path (6,0) coordinate (i);
\path (9,0) coordinate (m);
\path (6,1.5) coordinate (h);
\path (9,1.5) coordinate (l);
\draw (a) -- (b);
\draw (a) -- (c);
\draw (a) -- (e);
\draw (b) -- (c);
\draw (b) -- (f);
\draw (c) -- (d);
\draw (d) -- (e);
\draw (d) -- (f);
\draw (e) -- (f);
\draw (g) -- (h);
\draw (g) -- (i);
\draw (g) -- (j);
\draw (h) -- (j);
\draw (h) -- (l);
\draw (i) -- (j);
\draw (i) -- (m);
\draw (k) -- (l);
\draw (k) -- (m);
\draw (k) -- (n);
\draw (l) -- (n);
\draw (m) -- (n);
\draw (a) [fill=black] circle (\vr);
\draw (b) [fill=black] circle (\vr);
\draw (c) [fill=black] circle (\vr);
\draw (d) [fill=white] circle (\vr);
\draw (e) [fill=white] circle (\vr);
\draw (f) [fill=white] circle (\vr);
\draw (g) [fill=black] circle (\vr);
\draw (h) [fill=black] circle (\vr);
\draw (i) [fill=white] circle (\vr);
\draw (j) [fill=white] circle (\vr);
\draw (k) [fill=white] circle (\vr);
\draw (l) [fill=black] circle (\vr);
\draw (m) [fill=white] circle (\vr);
\draw (n) [fill=black] circle (\vr);
\draw (1.5,-0.75) node {(a) {\small $C_3 \cp K_2$}};
\draw (7.5,-0.75) node {(b) {\small $N_2$}};
\end{tikzpicture}
\end{center}
\vskip -0.65cm
\caption{The prism $C_3 \cp K_2$ and the diamond-necklace $N_2$.} \label{f:prism}
\end{figure}
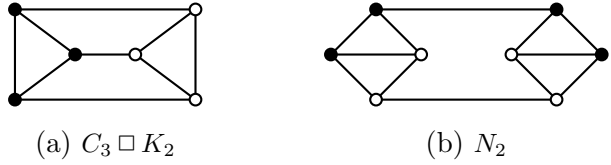

\medskip
\proof  We proceed by induction on the order~$n \ge 6$ of a connected, claw-free, cubic graph $G$. If $n = 6$, then $G$ is the prism $C_3 \cp K_2$ and $F_t(G) = 3 = n/2$ (the three darkened vertices shown in Figure~\ref{f:prism}(a) form a minimum TF-set in $C_3 \cp K_2$). If $n = 8$, then $G$ is the diamond-necklace $N_2$ and $F_t(G) = 4 = n/2$ (the four darkened vertices shown in Figure~\ref{f:prism}(b) form a minimum TF-set in $N_2$). This establishes the base cases. Let $n \ge 10$ and assume that if $G'$ is a connected, claw-free, cubic graph of order~$n'$, where $6 \le n' < n$, then $F_t(G') \le n'/2$ with equality if and only if $G' \in \cN_\cub$ or $G'$ is the prism $C_3 \cp K_2$. Let $G$ be a connected, claw-free, cubic graph of order~$n$.

\medskip

If $G \in \cN_\cub$, then by Lemma~\ref{lem1}, $F_t(G) = n/2$. Hence, we may assume that $G \notin \cN_\cub$, for otherwise the desired result follows. Therefore at least one unit in the $\Delta$-D-partition of $G$ is a triangle-unit. Since every triangle-unit is joined by three edges to vertices from other units, while every diamond-unit is joined by two edges to vertices from other units, we note that there are at least two triangle-units in our $\Delta$-D-partition. If $G$ has no diamond-unit, then the graph $G$ has a spanning $2$-factor consisting of triangles, and so by Theorem~\ref{thm1} $F_t(G) \le n/2$, with equality if and only if $G$ is the prism $C_3 \cp K_2$. Hence, we may assume that $G$ contains at least one diamond-unit, for otherwise the desired result follows.
Let $D$ be an arbitrary diamond-unit in the $\Delta$-D-partition of $G$, where $V(D) = \{a,b,c,d\}$ and where $ab$ is the missing edge in $D$. Let $e$ be the neighbor of $a$ not in $D$, and let $f$ be the neighbor of $b$ not in $D$. Since $G$ is claw-free, we note that $e \ne f$. We proceed further with the following claim.

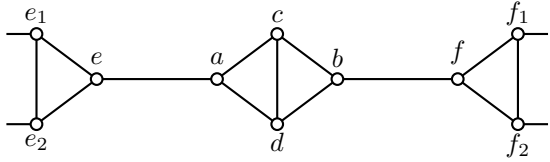
\begin{figure}[htb]
\begin{center}
\begin{tikzpicture}[scale=.8,style=thick,x=1cm,y=1cm]
\def\vr{2.75pt} 
\path (-0.5,0) coordinate (e2p);
\path (0,0) coordinate (e2);
\path (-0.5,1.5) coordinate (e1p);
\path (0,1.5) coordinate (e1);
\path (1,0.75) coordinate (e);
\path (3,0.75) coordinate (a);
\path (4,0) coordinate (d);
\path (4,1.5) coordinate (c);
\path (5,0.75) coordinate (b);
\path (7,0.75) coordinate (f);
\path (8.5,1.5) coordinate (f1p);
\path (8,1.5) coordinate (f1);
\path (8.5,0) coordinate (f2p);
\path (8,0) coordinate (f2);
\draw (e1) -- (e1p);
\draw (e2) -- (e2p);
\draw (f1) -- (f1p);
\draw (f2) -- (f2p);
\draw (e) -- (e1);
\draw (e) -- (e2);
\draw (e1) -- (e2);
\draw (f1) -- (f2);
\draw (e) -- (a);
\draw (f) -- (f1);
\draw (f) -- (f2);
\draw (f) -- (b);
\draw (a) -- (c);
\draw (a) -- (d);
\draw (b) -- (c);
\draw (b) -- (d);
\draw (c) -- (d);
\draw (a) [fill=white] circle (\vr);
\draw (b) [fill=white] circle (\vr);
\draw (c) [fill=white] circle (\vr);
\draw (d) [fill=white] circle (\vr);
\draw (e) [fill=white] circle (\vr);
\draw (e1) [fill=white] circle (\vr);
\draw (e2) [fill=white] circle (\vr);
\draw (f) [fill=white] circle (\vr);
\draw (f1) [fill=white] circle (\vr);
\draw (f2) [fill=white] circle (\vr);
\draw (0,-0.3) node {{\small $e_2$}};
\draw (0,1.8) node {{\small $e_1$}};
\draw (1,1.1) node {{\small $e$}};
\draw (3,1.1) node {{\small $a$}};
\draw (4,1.8) node {{\small $c$}};
\draw (4,-0.3) node {{\small $d$}};
\draw (5,1.1) node {{\small $b$}};
\draw (7,1.2) node {{\small $f$}};
\draw (8,1.85) node {{\small $f_1$}};
\draw (8,-0.35) node {{\small $f_2$}};
\end{tikzpicture}
\end{center}
\vskip -0.55cm
\caption{A subgraph of $G$.} \label{f:ClaimA}
\end{figure}

\begin{unnumbered}{Claim~A}
If $e$ and $f$ are not adjacent, then $F_t(G) < \frac{1}{2}n$.
\end{unnumbered}
\proof Suppose that $e$ and $f$ are not adjacent. Let $e_1$ and $e_2$ be the two neighbors of $e$ different from $a$, and let $f_1$ and $f_2$ be the two neighbors of $f$ different from $b$. By the claw-freeness of $G$, we note that the three vertices $e, e_1, e_2$ induce a triangle, say $T_e$, in $G$ and the three vertices $f, f_1, f_2$ induce a triangle, say $T_f$, in $G$. We note that $T_e$ and $T_f$ have no vertex in common. Further, $T_e$ (respectively, $T_f$) is either a triangle-unit or forms part of a diamond-unit. The resulting
subgraph of $G$ is illustrated in Figure~\ref{f:ClaimA}.

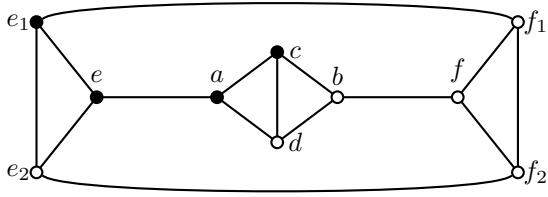
\begin{figure}[htb]
\begin{center}
\begin{tikzpicture}[scale=.8,style=thick,x=1cm,y=1cm]
\def\vr{2.75pt} 
\path (0,-0.5) coordinate (e2);
\path (0.1,-0.55) coordinate (e2p);
\path (0,2) coordinate (e1);
\path (0.1,2.05) coordinate (e1p);
\path (1,0.75) coordinate (e);
\path (3,0.75) coordinate (a);
\path (4,0) coordinate (d);
\path (4,1.5) coordinate (c);
\path (5,0.75) coordinate (b);
\path (7,0.75) coordinate (f);
\path (7.9,2.05) coordinate (f1p);
\path (8,2) coordinate (f1);
\path (7.9,-0.55) coordinate (f2p);
\path (8,-0.5) coordinate (f2);
%
\draw (e) -- (e1);
\draw (e) -- (e2);
\draw (e1) -- (e2);
\draw (f1) -- (f2);
\draw (e) -- (a);
\draw (f) -- (f1);
\draw (f) -- (f2);
\draw (f) -- (b);
\draw (a) -- (c);
\draw (a) -- (d);
\draw (b) -- (c);
\draw (b) -- (d);
\draw (c) -- (d);
%
\draw (a) [fill=black] circle (\vr);
\draw (b) [fill=white] circle (\vr);
\draw (c) [fill=black] circle (\vr);
\draw (d) [fill=white] circle (\vr);
\draw (e) [fill=black] circle (\vr);
\draw (e1) [fill=black] circle (\vr);
\draw (e2) [fill=white] circle (\vr);
\draw (f) [fill=white] circle (\vr);
\draw (f1) [fill=white] circle (\vr);
\draw (f2) [fill=white] circle (\vr);
\draw (-0.3,-0.5) node {{\small $e_2$}};
\draw (-0.3,2) node {{\small $e_1$}};
\draw (1,1.1) node {{\small $e$}};
\draw (3,1.1) node {{\small $a$}};
\draw (4.3,1.5) node {{\small $c$}};
\draw (4.3,0) node {{\small $d$}};
\draw (5,1.1) node {{\small $b$}};
\draw (7,1.2) node {{\small $f$}};
\draw (8.3,2) node {{\small $f_1$}};
\draw (8.3,-0.5) node {{\small $f_2$}};
\draw (e1p) to[out=45,in=135, distance=0.5cm] (f1p);
\draw (e2p) to[out=-45,in=-135, distance=0.5cm] (f2p);
\end{tikzpicture}
\end{center}
\vskip -0.55cm
\caption{The graph $G$.} \label{f:ClaimAa}
\end{figure}

Let $G'$ be the graph obtained from $G$ by deleting the vertices in $V(D)$ and their incident edges, and then adding the edge $ef$; that is, $G' = (G - V(D)) \cup \{ef\}$. We note that $G'$ is a connected, claw-free, cubic graph of order~$n'$, where $6 \le n' < n$. If $G' \in \cN_\cub$, then $G \in \cN_\cub$, contradicting our earlier assumption that $G \notin \cN_\cub$. Thus, $G' \notin \cN_\cub$. If $G'$ is the prism $C_3 \cp K_2$, then renaming vertices if necessary, we may assume that $e_i$ and $f_i$ are adjacent for $i \in [2]$, implying that $G$ is the graph of order~$n = 10$ shown in Figure~\ref{f:ClaimAa}. The set $\{a,c,e,e_1\}$, for example, illustrated in Figure~\ref{f:ClaimAa} is a TF-set in $G$, implying that $F_t(G) \le 4 < n/2$. Hence, we may assume that $G'$ is not the prism $C_3 \cp K_2$, for otherwise the desired result follows.

Applying the inductive hypothesis to the graph $G'$, we have $F_t(G') < n'/2$, noting that $G' \notin \cN_\cub$ and $G'$ is not the prism $C_3 \cp K_2$. We show next that $F_t(G) \le F_t(G') + 2$. Let $S'$ be a minimum TF-set in $G'$, and so $F_t(G') = |S'|$.

Suppose that $f \in S'$ but neither $f_1$ nor $f_2$ belong to $S'$. Since $G'[S']$ contains no isolated vertex, we note that in this case $e \in S'$. Thus, $S = (S' \setminus \{f\}) \cup \{a,c\}$ is a TF-set on $G$, since as the first three vertices played in the forcing process in $G$ starting with the set $S$ we play first the vertex $a$ which forces $d$ to be colored, second the vertex $d$ which forces $b$ to be colored, and third the vertex $b$ which forces $f$ to be colored, and then proceed with the forcing process in $G'$ starting with the set $S'$ to color all remaining vertices of $G$, implying that $F_t(G) \le |S| = |S'| + 1 = F_t(G') + 1$. Hence, we may assume that if $f \in S'$, then at least one of $f_1$ and $f_2$ belong to $S'$. Analogously, we may assume that if $e \in S'$, then at least one of $e_1$ and $e_2$ belong to $S'$.

Suppose that $e \in S'$. By our earlier assumptions, if $f \in S'$, then at least one of $f_1$ and $f_2$ belong to $S'$. Thus, the set $S = S' \cup \{a,c\}$ is a TF-set on $G$, since as the first two vertices played in the forcing process in $G$ starting with the set $S$ we play first the vertex $a$ which forces $d$ to be colored, and second the vertex $d$ which forces $b$ to be colored, and then proceed with the forcing process in $G'$ starting with the set $S'$ to color all remaining vertices of $G$, implying that $F_t(G) \le |S| = |S'| + 2 = F_t(G') + 2$. Analogously, if $f \in S'$, then $F_t(G) \le F_t(G') + 2$.

Suppose that neither $e$ nor $f$ belong to $S'$. In the forcing process in $G'$ starting with the set $S'$, we may assume renaming vertices if necessary that the vertex $e$ is colored before the vertex $f$. With this assumption, the set $S = S' \cup \{a,c\}$ is once again a TF-set on $G$, since we proceed with the forcing process in $G'$ starting with the set $S'$ until the vertex $e$ is colored, and then play the vertex $a$ which forces $d$ to be colored, followed by the vertex $d$ which forces $b$ to be colored, and thereafter continue with the original forcing process in $G'$ except that if the vertex $e$ is played in the original forcing process we instead play the vertex $b$. This implies that $F_t(G) \le |S| = |S'| + 2 = F_t(G') + 2$.

In all the above cases, we have shown that $F_t(G) \le F_t(G') + 2$. Hence since $F_t(G') < n'/2$ and $n' = n - 4$, this implies that $F_t(G) < n/2$. This completes the proof of Claim~A.~\smallqed

\begin{figure}[htb]
\begin{center}
\begin{tikzpicture}[scale=.8,style=thick,x=1cm,y=1cm]
\def\vr{2.75pt} 
\path (0.25,0.75) coordinate (c);
\path (1,1.5) coordinate (a);
\path (1,0) coordinate (b);
\path (1.75,0.75) coordinate (d);
\path (3,1.5) coordinate (e);
\path (3,0) coordinate (f);
\path (4,0.75) coordinate (g);
\path (5,0.75) coordinate (h);
\path (6,1.5) coordinate (i);
\path (6,0) coordinate (j);
\path (7,1.5) coordinate (k);
\path (7,0) coordinate (l);
\path (7.5,1.75) coordinate (k1);
\path (7.5,1.25) coordinate (k2);
\path (7.5,0.25) coordinate (l1);
\path (7.5,-0.25) coordinate (l2);
\draw (a) -- (e);
\draw (a) -- (c);
\draw (a) -- (d);
\draw (b) -- (c);
\draw (b) -- (d);
\draw (c) -- (d);
\draw (b) -- (f);
\draw (e) -- (f);
\draw (g) -- (e);
\draw (g) -- (f);
\draw (g) -- (h);
\draw (h) -- (i);
\draw (h) -- (j);
\draw (i) -- (k);
\draw (i) -- (j);
\draw (j) -- (l);
\draw (k) -- (k1);
\draw (k) -- (k2);
\draw (l) -- (l1);
\draw (l) -- (l2);
\draw (a) [fill=white] circle (\vr);
\draw (b) [fill=white] circle (\vr);
\draw (c) [fill=white] circle (\vr);
\draw (d) [fill=white] circle (\vr);
\draw (e) [fill=white] circle (\vr);
\draw (f) [fill=white] circle (\vr);
\draw (g) [fill=white] circle (\vr);
\draw (h) [fill=white] circle (\vr);
\draw (i) [fill=white] circle (\vr);
\draw (j) [fill=white] circle (\vr);
\draw (k) [fill=white] circle (\vr);
\draw (l) [fill=white] circle (\vr);
\draw (-0.05,0.75) node {{\small $c$}};
\draw (1,-0.35) node {{\small $b$}};
\draw (1,1.8) node {{\small $a$}};
\draw (2.05,0.75) node {{\small $d$}};
\draw (3,-0.35) node {{\small $f$}};
\draw (3,1.8) node {{\small $e$}};
\draw (4,0.4) node {{\small $g$}};
\draw (5,0.4) node {{\small $h$}};
\draw (6,-0.35) node {{\small $j$}};
\draw (6,1.85) node {{\small $i$}};
\draw (7,-0.35) node {{\small $\ell$}};
\draw (7,1.85) node {{\small $k$}};
\end{tikzpicture}
\end{center}
\vskip -0.55cm
\caption{A subgraph of $G$.} \label{f:figC}
\end{figure}
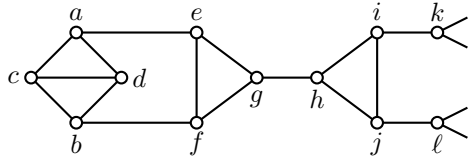

\medskip
By Claim~A, we may assume that $e$ and $f$ are adjacent, for otherwise $F_t(G) < n/2$, as desired. Thus, $e$ and $f$ belong to a common triangle-unit, say $T$. Let $g$ be the third vertex in the triangle-unit $T$, and let $h$ be the neighbor of $g$ not in $T$. Further, let $i$ and $j$ be the two vertices that belong to the triangle containing~$h$. If $h$ belongs to a diamond-unit, then choosing this diamond-unit as our initial diamond-unit $D$, we are back in Claim~A, implying that $F_t(G) < n/2$. Hence, we may assume that $h$ belongs to a triangle-unit in the $\Delta$-D-partition. Let $i$ and $j$ be the names of the vertices in this triangle-unit different from $h$. Further, let $k$ and $\ell$ denote the neighbors of $i$ and $j$, respectively, not in this triangle-unit. By assumption, $h$ does not belong to a diamond-unit, and so $k \ne \ell$. The resulting subgraph of $G$ is illustrated in Figure~\ref{f:figC}, where
possibly $k$ and $\ell$ are adjacent vertices.

\begin{figure}[htb]
\begin{center}
\begin{tikzpicture}[scale=.8,style=thick,x=1cm,y=1cm]
\def\vr{2.75pt} 
\path (0.25,0.75) coordinate (c);
\path (1,1.5) coordinate (a);
\path (1,0) coordinate (b);
\path (1.75,0.75) coordinate (d);
\path (3,1.5) coordinate (e);
\path (3,0) coordinate (f);
\path (4,0.75) coordinate (g);
\path (5,0.75) coordinate (h);
\path (6,1.5) coordinate (i);
\path (6,0) coordinate (j);
\path (8,1.5) coordinate (k);
\path (8,0) coordinate (l);
\path (7.25,0.75) coordinate (p);
\path (8.75,0.75) coordinate (m);
\draw (a) -- (e);
\draw (a) -- (c);
\draw (a) -- (d);
\draw (b) -- (c);
\draw (b) -- (d);
\draw (c) -- (d);
\draw (b) -- (f);
\draw (e) -- (f);
\draw (g) -- (e);
\draw (g) -- (f);
\draw (g) -- (h);
\draw (h) -- (i);
\draw (h) -- (j);
\draw (i) -- (k);
\draw (i) -- (j);
\draw (j) -- (l);
\draw (k) -- (m);
\draw (k) -- (p);
\draw (l) -- (m);
\draw (l) -- (p);
\draw (m) -- (p);
\draw (a) [fill=black] circle (\vr);
\draw (b) [fill=white] circle (\vr);
\draw (c) [fill=black] circle (\vr);
\draw (d) [fill=white] circle (\vr);
\draw (e) [fill=black] circle (\vr);
\draw (f) [fill=white] circle (\vr);
\draw (g) [fill=white] circle (\vr);
\draw (h) [fill=white] circle (\vr);
\draw (i) [fill=black] circle (\vr);
\draw (j) [fill=white] circle (\vr);
\draw (k) [fill=black] circle (\vr);
\draw (l) [fill=white] circle (\vr);
\draw (m) [fill=black] circle (\vr);
\draw (p) [fill=white] circle (\vr);
\draw (-0.05,0.75) node {{\small $c$}};
\draw (1,-0.35) node {{\small $b$}};
\draw (1,1.8) node {{\small $a$}};
\draw (2.05,0.75) node {{\small $d$}};
\draw (3,-0.35) node {{\small $f$}};
\draw (3,1.8) node {{\small $e$}};
\draw (4,0.4) node {{\small $g$}};
\draw (5,0.4) node {{\small $h$}};
\draw (6,-0.35) node {{\small $j$}};
\draw (6,1.85) node {{\small $i$}};
\draw (8,-0.35) node {{\small $\ell$}};
\draw (8,1.85) node {{\small $k$}};
\draw (6.95,0.7) node {{\small $p$}};
\draw (9.1,0.75) node {{\small $m$}};
\end{tikzpicture}
\end{center}
\vskip -0.55cm
\caption{The graph $G$.} \label{f:Gdeter3}
\end{figure}
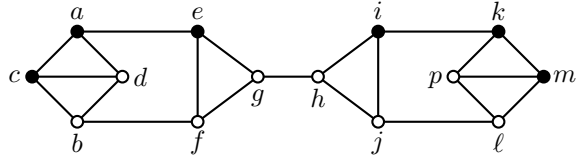

We now consider the connected, claw-free, cubic graph $G'$ obtained from $G$ by removing the vertices of the set $\{e,f,g,h,i,j\}$ and adding the edges $ak$ and $b\ell$. Let $n'$ be the order of $G'$, and so $n' = n - 6$. Suppose that the vertex $k$ belongs to a diamond-unit $D^*$ in $G$ (and therefore in $G'$). If this diamond-unit does not contain the vertex~$\ell$, then choosing this diamond-unit as our initial diamond-unit $D$, we are back in Claim~A, implying that $F_t(G) < n/2$. Hence, we may assume that the diamond-unit $D^*$ contains the vertex $\ell$, implying that the graph $G$ is the graph of order~$n = 14$ shown in Figure~\ref{f:Gdeter3}, where $V(D^*) = \{k,\ell,m,p\}$. The set $\{a,c,e,i,k,m\}$, for example, illustrated in Figure~\ref{f:Gdeter3} is a TF-set in $G$, implying that $F_t(G) \le 6 < n/2$. Hence, we may assume that the vertex $k$ belongs to a triangle-unit in $G$ (and therefore in $G'$). Thus, $G' \notin \cN_\cub$. Since $G'$ contains the diamond-unit $D$, we note that $G'$ is not the prism $C_3 \cp K_2$. Applying the inductive hypothesis to the graph $G'$, we therefore have $F_t(G') < n'/2$.

We show next that $F_t(G) \le F_t(G') + 3$. Let $S'$ be a minimum TF-set in $G'$, and so $F_t(G') = |S'|$. We note that every diamond-unit in $G'$ contains at least two vertices of $S'$. In particular, the diamond-unit $D$ contains at least two vertices of $S'$. We now consider the set $S$ obtained from $S'$ by removing the vertices of $D$ in $S'$ and adding to it the vertices in the set $\{a,c,e,i,j\}$; that is, $S = (S' \setminus V(D)) \cup \{a,c,e,i,j\}$. The set $S$ is a TF-set of $G$, since as the first five vertices played in the forcing process starting with the set $S$ we play first the vertex $a$ which forces $d$ to be colored, second we play the vertex $d$ which forces $b$ to be colored, third we play the vertex $b$ which forces $f$ to be colored, fourth we play the vertex $f$ which forces $g$ to be colored, fifth we play the vertex $g$ which forces $h$ to be colored, and thereafter we follow the forcing process starting with the set $S'$ in $G'$ to color all remaining vertices (where we replace the vertex $a$ with the vertex $i$ if $a$ is played to color $k$ in the original forcing process in $G'$ and we replace the vertex $b$ with the vertex $j$ if $b$ is played to color $\ell$ in the original forcing process in $G'$). Hence, $F_t(G) \le |S| \le |S'| - 3 < n'/2 - 3 = n/2$. This completes the proof of Theorem~\ref{thm2}.~\qed

\subsection{Proof of Theorem~\ref{thm3}}
\label{S:proof3}

In this section, we present a proof of Theorem~\ref{thm3}. Recall its statement.

\noindent \textbf{Theorem~\ref{thm3}}. \emph{If $G \ne K_4$ is a connected, claw-free, cubic graph of order $n$, then $F(G) \le \frac{1}{2}n$ with equality if and only if $G$ is the diamond-necklace $N_2$ or the prism $C_3 \cp K_2$.
}

\medskip
\proof  Let $G \ne K_4$ be a connected, claw-free, cubic graph of order $n$. By Theorem~\ref{thm2} and Observation~\ref{ob1}, $F(G) \le F_t(G) \le n/2$. Further if $F(G) = n/2$, then $F_t(G) = n/2$, implying by Theorem~\ref{thm2} that $G \in \cN_\cub$ or $G$ is the prism $C_3 \cp K_2$. Suppose that $G \in \cN_\cub$ has order~$n$, and so $n = 4k$ for some $k \ge 2$. If $n \ge 12$, then $F(G) = n/4 + 2 < n/2$, a contradiction. Hence, $n = 8$, and so, $G = N_2$. If $G = C_3 \cp K_2$, then $n = 6$ and $F(G) = 3 = n/2$.~\qed

\subsection{Proof of Theorem~\ref{thm4}}
\label{S:proof4}

In this section, we present a proof of Theorem~\ref{thm4}. Recall its statement.

\noindent \textbf{Theorem~\ref{thm4}}. \emph{If $G$ is a connected, claw-free, cubic graph of order $n$, then
\[
\frac{F_t(G)}{F(G)} \le 2,
\]
and this bound is asymptotically best possible. 
}

\medskip
\proof  By Observation~\ref{ob1}, $F_t(G)/F(G) \le 2$. To show that this bound is asymptotically best possible, let $\epsilon > 0$ be an arbitrary real number and consider the graph $G = N_k \in \cN_\cub$ where $k > \frac{4}{\epsilon} - 2$. By Lemma~\ref{lem1},
\[
\frac{F_t(G)}{F(G)} = \frac{2k}{k+2} = 2 - \frac{4}{k+2} > 2 - \epsilon,
\]
implying that the ratio $F_t(G)/F(G)$ can be arbitrarily close to~$2$.~\qed

\end{document}